.
\font\sets=msbm10.
\font\script=eusm10.
\font\stampatello=cmcsc10.
.
\def\1{{\bf 1}}
\def\defineq{\buildrel{def}\over{=}}
\def\definiz{\buildrel{def}\over{\Longleftrightarrow}}

\def\C{\hbox{\sets C}}
\def\N{\hbox{\sets N}}

\def\Z{\hbox{\sets Z}}

\def\0{{\bf 0}}

\def\Carmichael{{\rm Car}}
\def\Wintner{{\rm Win}}

\def\square{\hbox{\vrule\vbox{\hrule\phantom{s}\hrule}\vrule}}

\par
\centerline{\bf A smooth shift approach for a Ramanujan expansion}
\bigskip
\bigskip
\bigskip
\centerline{Giovanni Coppola}
\bigskip
\bigskip
\rightline{\it to oncoming smooth-numbers aficionados} 
\bigskip
\bigskip
\bigskip
\par
\noindent
{\bf Abstract}. All arithmetical functions $F$ satisfying Ramanujan Conjecture, i.e., $F(n)\ll_{\varepsilon}n^{\varepsilon}$, and with $Q-$smooth divisors, i.e., with Eratosthenes transform $F'\defineq F\ast \mu$ supported in $Q-$smooth numbers, have a kind of {\it unique} Ramanujan expansion; also, these Ramanujan coefficients decay very well to $0$ and have two explicit expressions (in the style of Carmichael and Wintner). This general result, then, is applied to the {\it shift-}Ramanujan expansions, i.e., the expansions for correlations with respect to the shift, whence the title. 
\bigskip
\bigskip
\bigskip
\par
\noindent {\bf 1. Introduction, statements and proofs of the results.}
\smallskip
\par
\noindent
In the following, we fix $Q\in \N$ and indicate the set of $Q-smooth$ (positive) integers writing 
$$
(Q)\defineq \{ n\in \N : n=1\enspace \hbox{\rm or}\enspace p|n\;\Rightarrow\;p\le Q\}
$$
\par
\noindent
(now on $p$ denotes a prime, eventually with subscripts) and writing (as usual $(a,b)\defineq g.c.d.(a,b)$ now on) 
$$
)Q(\,\defineq \{ n\in \N : (n,\prod_{p\le Q}p)=1\}
$$
\par
\noindent
the set of $Q-sifted$ (positive) integers. See that $(Q)\cap\,)Q(=\{1\}$ and $n\in(Q)$, $m\in\,)Q($ implies $(n,m)=1$. 
\par
\noindent
We need to define the $Q-smooth$ {\it restriction} of any $F:\N \rightarrow \C$ as 
$$
F_{(Q)}(n)\defineq \sum_{{d|n}\atop{d\in(Q)}}F'(d), 
\enspace \forall n\in \N, 
$$
\par
\noindent
where as usual $F'$ is the {\it Eratosthenes transform} [W] of $F$, namely \enspace $F'\defineq F\ast \mu$. (See [T] for $\ast$, Dirichlet product, and $\mu$, M\"{o}bius function.) 
\par
Notice, in passing, that the Eratosthenes transform of our $F_{(Q)}$, namely $(F_{(Q)})'$, thanks to 
$$
F_{(Q)}(n)=\sum_{d|n}F'(d)\1_{(Q)}(d), 
\enspace \forall n\in \N, 
$$
\par
\noindent
is nothing else than $F'\cdot \1_{(Q)}$, with $\1_{\hbox{\script A}}$ the characteristic function of the set $\hbox{\script A}$. (Here {\script A}$=\{Q-$smooth n.s.$\}$) 
\par
See the similarity of notation with $F_Q$, which is the $Q-$truncation of our $F$, namely we truncate its divisors after $Q$, i.e., the Eratosthenes transform, now, has support in $\{1,\ldots,Q\}$ (compare [C2], [CMS] and [CM]). Our $Q-$smooth restriction has an infinity of divisors, while of course $F_Q$ has only at most $Q$ of them! 
\par
While the $Q-$truncations (i.e., {\it truncated divisor sums}) are strictly connected to finite Ramanujan expansions (see section 5 of [C2] and compare [CMS], [CM]), here the $Q-$restrictions (i.e., {\it restricted divisor sums}) are linked, see $(RE)$ in next Theorem 1, to infinite, pointwise converging Ramanujan expansions ! 
\par
For $F:\N \rightarrow \C$ we define [C2] {\it Carmichael's coefficients} (provided following limits exist) and {\it Wintner's coefficients} (if following series converge), with $\varphi(q)$ the {\it Euler function} and $c_q(n)$ the {\it Ramanujan sum} [R],[M]: 
$$
\Carmichael_q(F)\defineq {1\over {\varphi(q)}}\lim_{x\to\infty}{1\over x}\sum_{n\le x}F(n)c_q(n), 
\enspace
\forall q\in \N, 
\qquad
\Wintner_q(F)\defineq \sum_{{{d=1}\atop {d\equiv 0(\!\!\bmod q)}}}^{\infty}{{F'(d)}\over d}, 
\enspace
\forall q\in \N. 
$$				
\par
The main limit of our Theorem 1 [CM] is that we need a hypothesis (we may choose among three) for the finite Ramanujan expansion of (suitable) shifted convolution sums (SCS); now, we don't need it, simply considering not the original SCS, but restricting its divisors, as above for $F$, and then it can be expanded into, say, the most expected but the smoothest Ramanujan expansion; most expected, as the coefficients are nothing else than Carmichael's \& Wintner's at the same time (like we expect from Wintner-Delange Formula, see Theorem 2 in [C2]) and smoothest, as the coefficients are so smooth that they satisfy what we call in [C2] the \lq \lq Dual Delange Assumption\rq \rq, which also guarantees the uniqueness of these Carmichael-Wintner, say, coefficients as the {\it unique} Ramanujan coefficients. (See Theor.1 for general $F$ and Corollary 1 for SCS.)
\bigskip
\par
We wish to prove, before general results in our Theorem 1, a Proposition that regards suitable SCS, also called {\it correlations}, that will be applied in our Corollary 1 (for correlations). 
\par
Its, say, Basic Hypothesis needs two definitions.
\par
We call a general arithmetic function $g:\N \rightarrow \C$ \lq \lq {\it of range } $Q$\rq \rq, by definition, when $g$ may be expressed through its Eratosthenes transform $g'$, for a fixed $Q\in \N$, as: (a truncated divisor sum!) 
$$
g(m)\defineq \sum_{d|n,d\le Q}g'(d),
\enspace \forall m\in \N. 
$$
\par
\noindent
(Compare, for a more rigorous definition, [C2].)
\par
Once given a {\it correlation} (or Shifted Convolution Sum, SCS) of fixed $f,g:\N \rightarrow \C$, i.e., 
$$
C_{f,g}(N,a)\defineq \sum_{n\le N}f(n)g(n+a),
\enspace \forall a\in \N,
$$
\par
\noindent
where the \lq \lq length\rq \rq, $N\in \N$, is fixed and the \lq \lq {\it shift}\rq \rq, $a\in \N$, is our variable (so that the Eratosthenes transform of $C_{f,g}(N,a)$ is $C'_{f,g}(N,d)\defineq \sum_{t|d}C_{f,g}(N,t)\mu(d/t)$, see  Corollary 1), we say that
\par
$C_{f,g}(N,a)$ is {\it fair} $\definiz$ dependence on the shift $a$ is only inside $g$'s argument $(n+a)$ 
\smallskip
\par
\noindent
(i.e. nor dependence on $a$ inside $f,g$, neither in their supports; esp., $f_H$'s correlation is not fair: [CM] end) 
\medskip
\par
We prove very quickly a property of correlations, in \lq \lq Basic Hypothesis\rq \rq, i.e., the two hypotheses of Theorem 1 [CM]; in fact, following Proposition is already \lq \lq implicit\rq \rq, in [CM] Theorem 1 Proof.
\smallskip
\par
\noindent {\bf Proposition 1.}
\enspace {\it Let } $f,g:\N \rightarrow \C$ {\it be such that }
$$
g\enspace \hbox{\it is } \hbox{\stampatello of range } Q\le N
\enspace \hbox{\it and }\enspace C_{f,g}(N,a) \enspace \hbox{\it is} \enspace \hbox{\stampatello fair}. 
\leqno{(BH)}
$$  
\par
\noindent
{\it Then}
\smallskip
\item{$(i)$} $C_{f,g}(N,a)={\displaystyle \sum_{q\le Q}\widehat{g}(q)\sum_{n\le N} }f(n)c_{q}(n+a)$, $\forall a\in \N$, 
\enspace {\it where } $\widehat{g}(q){\displaystyle \defineq \sum_{{d\le Q}\atop {d\equiv 0(\!\!\bmod q)}}{{g'(d)}\over d} }$, $\forall q\in \N$; 
\smallskip
\item{$(ii)$} $C_{f,g}(N,a)$ {\it is, with respect to } $a\in \N$, {\it periodic, whence bounded}; 
\smallskip
\item{$(iii)$} $C_{f,g}(N,a)$ {\it has coincident Carmichael and Wintner $\ell-$th coefficients:} ${\displaystyle {{\widehat{g}(\ell)}\over {\varphi(\ell)}}\sum_{n\le N}f(n)c_{\ell}(n) }$, $\forall \ell\in \N$. 
\smallskip
\par
\noindent {\bf Proof.} Here, $(i)$ follows from the $g$ finite Ramanujan expansion $g(n+a)=\sum_{q\le Q}\widehat{g}(q)c_q(n+a)$ of Ramanujan coefficients $\widehat{g}(q)$ as above, see [C2], beginning of section 5. 
\par
Then, from $(i)$, together with fairness we get, since each $c_q(n+a)$ is periodic modulo $q$, with respect to $a$, periodicity (with period dividing {\script Q}$\defineq $l.c.m.$(2,\ldots,Q)$, of course), w.r.t. $a$, whence $C_{f,g}(N,a)$ is bounded (w.r.t. $a$). 
\par
Finally, $(iii)$ follows from the Delange 1987 Theorem [De87], in the equivalent form, given as Theorem 9 in [C2]; in fact, our $C_{f,g}(N,a)$ is bounded, so bounded on average, as required by Th.9 first assumption and its second assumption is satisfied because $C_{f,g}(N,a)$ has all the Carmichael coefficients, since by $(i)$ and fairness we get (compare $(CC)$ in [C2]), as $\ell$-th Carmichael coefficient of $C_{f,g}(N,a)$,  
$$
{1\over {\varphi(\ell)}}\lim_x {1\over x}\sum_{a\le x}c_{\ell}(a)\sum_{q\le Q}\widehat{g}(q)\sum_{n\le N}f(n)c_{q}(n+a)
={1\over {\varphi(\ell)}}\sum_{q\le Q}\widehat{g}(q)\sum_{n\le N}f(n)\lim_x {1\over x}\sum_{a\le x}c_{\ell}(a)c_{q}(n+a), 
$$
\par
\noindent
whence, applying the {\it orthogonality of Ramanujan sums} (proved by Carmichael in 1932 [Ca], see Theorem 1 in [M]), we get $(iii)$.\hfill $\square$ 

\vfill
\eject

\par				
We come to our main result. Hereafter, we write $V$ to avoid confusion with $Q$ in Proposition 1, we assume $V\in \N$ and, as usual [Da], $\omega(d)$ will be the number of prime factors of $d$. 
\smallskip
\par
\noindent {\bf Theorem 1.}
\enspace {\it Let} $F:\N \rightarrow \C$ {\it satisfy Ramanujan Conjecture and fix an integer $V>1$. Then} 
\medskip
\item{\it (i)} $\Carmichael_{\ell}(F_{(V)})=\Wintner_{\ell}(F_{(V)})$, $\forall \ell\in \N$ \enspace {\it and in particular } $\Carmichael_{\ell}(F_{(V)})=\Wintner_{\ell}(F_{(V)})=0$, $\forall \ell\not\in (V)$; 
\medskip
\item{\it (ii)} ${\displaystyle \Carmichael_{\ell}(F_{(V)})=\prod_{p\le V}\left(1-{1\over p}\right){1\over {\varphi(\ell)}}\sum_{t\in(V)}{{F(t)c_{\ell}(t)}\over t} }$, $\forall \ell \in (V)$, {\it where } ${\displaystyle \sum_{t\in(V)}{{|F(t)c_{\ell}(t)|}\over t}<\infty }$, $\forall \ell \in \N$; 
\medskip
\item{\it (iii)} 
${\displaystyle F_{(V)}(a)=\sum_{\ell\in (V)}\left(\prod_{p\le V}\left(1-{1\over p}\right){1\over {\varphi(\ell)}}\sum_{t\in(V)}{{F(t)c_{\ell}(t)}\over t}\right)c_{\ell}(a)
=\sum_{\ell\in (V)}\left(\sum_{{d\in(V)}\atop{d\equiv 0(\!\!\bmod \ell)}}{{F'(d)}\over d}\right)c_{\ell}(a) }$, $\forall a\in \N$, 
\par
{\it whence, in particular,}
$$
F(a)=\sum_{\ell\in (V)}\left(\prod_{p\le V}\left(1-{1\over p}\right){1\over {\varphi(\ell)}}\sum_{t\in(V)}{{F(t)c_{\ell}(t)}\over t}\right)c_{\ell}(a)
=\sum_{\ell\in (V)}\left(\sum_{{d\in(V)}\atop{d\equiv 0(\!\!\bmod \ell)}}{{F'(d)}\over d}\right)c_{\ell}(a), 
\thinspace \forall a\in (V); 
\leqno{(RE)}
$$ 
\medskip
\item{\it (iv)} {\it the Ramanujan coefficients } $\widehat{F_{(V)}}(\ell)\defineq \Carmichael_{\ell}(F_{(V)})=\Wintner_{\ell}(F_{(V)})$ {\it satisfy } ${\displaystyle \sum_{\ell=1}^{\infty}2^{\omega(\ell)}|\widehat{F_{(V)}}(\ell)|<\infty }$; 
\medskip
\item{\it (v)} $F_{(V)}(a)=\sum_{\ell=1}^{\infty}R_{(V),F}(\ell)c_{\ell}(a)$, $\forall a\in \N$ {\it and $(iv)$ holds for } $R_{(V),F}(\ell)$ $\Rightarrow $ $R_{(V),F}(\ell)=\widehat{F_{(V)}}(\ell)$, $\forall \ell\in \N$.
\medskip
\par
\noindent {\bf Proof.} Before going on, recall from the definition that the Wintner coefficients of $F_{(V)}$ are 
$$
\Wintner_{\ell}(F_{(V)})=\sum_{{d\in(V)}\atop{d\equiv 0(\!\!\bmod \ell)}}{{F'(d)}\over d}, 
$$
\par
\noindent
in which, of course, the condition $\ell|d$ implies, since $d\in(V)$, that $\ell\in(V)$; otherwise, the coefficient vanishes. So we are left with the task to prove coincidence of Carmichael and Wintner $\ell-$th coefficients, for all $\ell\in\N$. 
We start proving this part of $(i)$. This can be done, proving the Delange Hypothesis 
$$
\sum_{d=1}^{\infty}{{2^{\omega(d)}|(F_{(V)})'(d)|}\over d}<\infty
\leqno{(DH)}
$$ 
\par
\noindent
because Delange 1976 Theorem [De] infers from $(DH)$ both the identity of Carmichael \& Wintner coefficients, i.e. $(i)$ (for what we saw above), and the convergence, of corresponding Ramanujan expansion; thus proving, after we prove $(ii)$, also $(iii)$. In order to prove $(DH)$ above, recall \enspace $(F_{(V)})'(d)=F'(d)\1_{(V)}(d)$ \enspace so (hereafter we use classic notation, Vinogradov's $\ll$ and Landau's $O$, like $\pi(V)\defineq |\{p\le V\}|$, see [Da])
$$
\sum_{d=1}^{\infty}{{2^{\omega(d)}|(F_{(V)})'(d)|}\over d}=\sum_{d\in(V)}{{2^{\omega(d)}|F'(d)|}\over d}
\ll 2^{\pi(V)}\sum_{d\in(V)}{{|F'(d)|}\over d}
\ll_{V,\varepsilon} \sum_{d\in(V)}d^{\varepsilon-1}
<\infty,
$$ 
\par
\noindent
where $F$ satisfying the Ramanujan Conjecture implies the same for $F'$, then we apply Lemma 3 (at next $\S2$). 
\smallskip
\par
We have proved both $(i)$ and $(iii)$, once we prove $(ii)$, too. 
\smallskip
\par
\noindent
For $(ii)$ we start proving the absolute convergence:
$$
\sum_{t\in(V)}{{|F(t)c_{\ell}(t)|}\over t}\ll_{\varepsilon} \sum_{t\in(V)}(\ell,t)t^{\varepsilon-1}
\ll_{\varepsilon} \sum_{{d\in(V)}\atop{d|\ell}}d\sum_{{t\in(V)}\atop{t\equiv 0(\!\!\bmod d)}}t^{\varepsilon-1}
\ll_{\varepsilon} \sum_{{d\in(V)}\atop{d|\ell}}d^{\varepsilon}\sum_{K\in(V)}K^{\varepsilon-1}
\ll_{V,\varepsilon,\ell}1, 
$$
\par				
\noindent
using the inequality $|c_{\ell}(t)|\le (\ell,t)$ (see Lemma A.1 in [CM]) and Lemma 3 at $\S2$. We have left to prove the formula for Carmichael $\ell-$th coefficient, once $\ell\in(V)$. Starting with the definition for these coefficients and adding M\"{o}bius switch, namely Lemma 1 in next $\S2$, 
$$
\Carmichael_{\ell}(F_{(V)})={1\over {\varphi(\ell)}}\lim_x {1\over x}\sum_{a\le x}c_{\ell}(a)\sum_{{d\in(V)}\atop{d|a}}F'(d)
={1\over {\varphi(\ell)}}\lim_x \sum_{{t\in(V)}\atop{t\le x}}F(t)\cdot {1\over x}\sum_{{K\in\,)V(}\atop{K\le {x\over t}}}c_{\ell}(tK), 
$$
\par
\noindent
true $\forall \ell\in \N$; however, assuming $\ell\in(V)$ now, we use the fact that $K\in\,)V($ to get $(\ell,K)=1$, whence $c_{\ell}(tK)=c_{\ell}(t)$, $\forall t\in \N$, getting from the count in Lemma 2 ($\S2$) 
$$
\Carmichael_{\ell}(F_{(V)})={1\over {\varphi(\ell)}}\lim_x \sum_{{t\in(V)}\atop{t\le x}}F(t)c_{\ell}(t)\cdot {1\over x}\sum_{{K\in\,)V(}\atop{K\le {x\over t}}}1
={1\over {\varphi(\ell)}}\lim_x \sum_{{t\in(V)}\atop{t\le x}}F(t)c_{\ell}(t)\left({1\over t}\prod_{p\le V}\left(1-{1\over p}\right)+O_V\left({1\over x}\right)\right), 
$$
\par
\noindent
in which, using the absolute convergence just proved, we have convergence of main term, i.e. 
$$
{1\over {\varphi(\ell)}}\lim_x \sum_{{t\in(V)}\atop{t\le x}}{{F(t)c_{\ell}(t)}\over t}\prod_{p\le V}\left(1-{1\over p}\right)=\prod_{p\le V}\left(1-{1\over p}\right){1\over {\varphi(\ell)}}\sum_{t\in(V)}{{F(t)c_{\ell}(t)}\over t}, 
$$
\par
\noindent
so, we're left with proving that remainders don't count, as next term is infinitesimal with $x\to \infty$ :  
$$
\sum_{{t\in(V)}\atop{t\le x}}F(t)c_{\ell}(t)O_V\left({1\over x}\right)\ll_{V,\varepsilon} x^{\varepsilon-1}\sum_{{t\in(V)}\atop{t\le x}}(\ell,t)
\ll_{V,\varepsilon} x^{\varepsilon-1}\sum_{d|\ell}d\sum_{{K\in(V)}\atop{K\le x/d}}1
\ll_{V,\varepsilon} x^{\varepsilon-1}\sum_{d|\ell}d\left({x\over d}\right)^{\varepsilon}
\ll_{V,\varepsilon,\ell} x^{2\varepsilon-1}, 
$$
\par
\noindent
applying, in penultimate step, the bound of Lemma 2 ($\S2$) and recalling $\varepsilon>0$ is small, finally proving $(ii)$.
\smallskip
\par
We come to proving $(iv)$, now. 
\smallskip
\par
\noindent
Since 
$$
\sum_{\ell=1}^{\infty}2^{\omega(\ell)}|\widehat{F_{(V)}}(\ell)|=\sum_{\ell\in(V)}2^{\omega(\ell)}\left|\sum_{{d\in(V)}\atop{d\equiv 0(\!\!\bmod \ell)}}{{F'(d)}\over d}\right|
\le 2^{\pi(V)}\sum_{\ell\in(V)}\sum_{{d\in(V)}\atop{d\equiv 0(\!\!\bmod \ell)}}{{|F'(d)|}\over d}
\ll_{V,\varepsilon} \sum_{\ell\in(V)}\ell^{\varepsilon-1}\sum_{K\in(V)}K^{\varepsilon-1}, 
$$
\par
\noindent
we prove, here, even more than $(iv)$, thanks to Lemma 3 (at $\S2$), again. 
\par
We need only to prove $(v)$, a kind of \lq \lq uniqueness\rq \rq, for the Ramanujan expansion we found (with Carmichael coefficients $=$ Wintner coefficients). It follows from Theorem 4 of [C2], recalling $(iv)$ is a kind of \lq \lq Dual Delange\rq \rq, as we call it in [C2], assumption.\hfill $\square$ 

\medskip

\par
Combining Theorem 1 and $(ii)$ of Proposition 1 we easily get the following (which we don't prove). 
\smallskip
\par
\noindent {\bf Corollary 1.}
\enspace {\it Given } $f:\N \rightarrow \C$ {\it and } $g:\N \rightarrow \C$, {\it satisfying the Basic Hypothesis $(BH)$ above, then, given any integer $V>1$, for } $G_{(V),f,g,N}(a)=G(a)\defineq \sum_{d|a,d\in(V)}C'_{f,g}(N,d)=\sum_{d|a}C'_{f,g}(N,d)\1_{(V)}(d)$, $\forall a\in\N$ {\it we have the following } {\stampatello Ramanujan expansion} 
$$
G(a)=\sum_{\ell\in (V)}\Big(\sum_{{d\in(V)}\atop{d\equiv 0(\!\!\bmod \ell)}}{{C'_{f,g}(N,d)}\over d}\Big)c_{\ell}(a) 
=\sum_{\ell\in (V)}\Big(\prod_{p\le V}\left(1-{1\over p}\right){1\over {\varphi(\ell)}}\sum_{t\in(V)}{{C_{f,g}(N,t)c_{\ell}(t)}\over t}\Big)c_{\ell}(a)
, \enspace \forall a\in \N, 
$$
\par
\noindent
{\it whence, in particular, } $\forall a\in (V)$, 
$$
C_{f,g}(a)=\sum_{\ell\in (V)}
\sum_{{d\in(V)}\atop{d\equiv 0(\!\!\bmod \ell)}}{{C'_{f,g}(N,d)}\over d}
c_{\ell}(a) 
=\sum_{\ell\in (V)}
\prod_{p\le V}\left(1-{1\over p}\right){1\over {\varphi(\ell)}}\sum_{t\in(V)}{{C_{f,g}(N,t)c_{\ell}(t)}\over t}
c_{\ell}(a). 
\leqno{(SCS)}
$$ 
\par				
\noindent
{\it The \lq \lq Dual Delange\rq \rq \thinspace assumption on the Ramanujan coefficients, say $\widehat{G}(\ell)$, of $G$ :}
$$
\sum_{\ell=1}^{\infty}2^{\omega(\ell)}|\widehat{G}(\ell)|<\infty
\leqno{(DD)}
$$
\par
\noindent
{\it holds for the coefficients above (which vanish outside of $V-$smooth numbers) and $\underline{\hbox{ONLY}}$ for them:
 IF}
$$
G(a)=\sum_{\ell=1}^{\infty}\widehat{G}(\ell)c_{\ell}(a), \forall a\in\N
\enspace \hbox{\stampatello and } (DD) \enspace \hbox{\stampatello holds } 
$$
\par
\noindent
{\it THEN}
$$
\widehat{G}(\ell)=\sum_{{d\in(V)}\atop{d\equiv 0(\!\!\bmod \ell)}}{{C'_{f,g}(N,d)}\over d}
=\prod_{p\le V}\left(1-{1\over p}\right){{\1_{(V)}(\ell)}\over {\varphi(\ell)}}\sum_{t\in(V)}{{C_{f,g}(N,t)c_{\ell}(t)}\over t}, 
\enspace \forall \ell\in \N.  
$$
\medskip
\par
\noindent {\bf Remark 1.} The same thesis comes from alternative hypothesis: $f$ \& $g$ satisfy Ramanujan Conjecture. 
\medskip
\par
\noindent {\bf Corollary 2.} {\it Given } $f,g:\N \rightarrow \C$ {\it satisfying } $(BH)$, {\it then we have, } $\forall V>1$ {\it integer}, 
$$
\sum_{{d\in(V)}\atop{d\equiv 0(\!\!\bmod \ell)}}{{C'_{f,g}(N,d)}\over d}
={{\widehat{g}(\ell)}\over {\varphi(\ell)}}\sum_{n\le N}f(n)c_{\ell}(n)
 -\sum_{{d\not\in(V)}\atop{d\equiv 0(\!\!\bmod \ell)}}{{C'_{f,g}(N,d)}\over d}, 
\enspace \forall \ell\in \N, 
$$
\par
\noindent
{\it which, in particular, gives}
$$
\sum_{{d\not\in(V)}\atop{d\equiv 0(\!\!\bmod \ell)}}{{C'_{f,g}(N,d)}\over d}={{\widehat{g}(\ell)}\over {\varphi(\ell)}}\sum_{n\le N}f(n)c_{\ell}(n), 
\enspace \forall \ell\not\in(V). 
$$
\smallskip
\par
\noindent {\bf Proof.} The application of $(BH)$ in Proposition 1 gives $(iii)$, implying 
$$
\sum_{d\equiv 0(\!\!\bmod \ell)}{{C'_{f,g}(N,d)}\over d} = {{\widehat{g}(\ell)}\over {\varphi(\ell)}}\sum_{n\le N}f(n)c_{\ell}(n),
\enspace \forall \ell \in \N. 
$$
\par
\noindent
Then, we may separate $d\in(V)$ and $d\not\in(V)$ series, thanks to absolute convergence in Wintner coefficients, say, with Eratosthenes transform restricted to $V-$smooth numbers, i.e., 
$$
\sum_{{d\in(V)}\atop{d\equiv 0(\!\!\bmod \ell)}}{{|C'_{f,g}(N,d)|}\over d}<\infty;
\leqno{(\ast)} 
$$
\par
\noindent
in fact : $C_{f,g}(N,a)$ bounded $\Rightarrow $ $C'_{f,g}(N,d)\ll_{N,Q,\varepsilon} d^{\varepsilon}$, whence 
$$
\sum_{{d\in(V)}\atop{d\equiv 0(\!\!\bmod \ell)}}{{|C'_{f,g}(N,d)|}\over d}
\ll_{N,Q,\varepsilon} \sum_{{d\in(V)}\atop{d\equiv 0(\!\!\bmod \ell)}}d^{\varepsilon-1}
\ll_{N,Q,\varepsilon} \ell^{\varepsilon-1}\sum_{K\in(V)}K^{\varepsilon-1}
\ll_{N,Q,V,\varepsilon} \ell^{\varepsilon-1}
\ll_{N,Q,V,\varepsilon,\ell} 1, 
$$
\par
\noindent
implying $(\ast)$ above, from Lemma 3 at next $\S2$.\hfill $\square$ 
\medskip
\par
\noindent {\bf Remark 2.} The alternative hypothesis, $f$ \& $g$ satisfy Ramanujan Conjecture, this time doesn't suffice (as we are using $(BH)$ for Wintner coefficients formula). 
\medskip

\vfill

We give the Lemmas used above, in next $\S2$; then, in $\S3$ a kind of new orthogonality relations for Ramanujan sums provide a new approach to Theorem 1, see Proposition 2. Our Conjectures (compare version two) are disproved in $\S4$; finally, $\S5$ gives  further remarks. 

\eject

\par				
\noindent {\bf 2. Lemmas.}
\smallskip
\par
\noindent
We give a page of Lemmas for our proofs. 

\bigskip

\par
First Lemma is \lq \lq M\"{o}bius switch\rq \rq.
\smallskip
\par
\noindent
{\bf Lemma 1.}
\enspace {\it For any $F:\N \rightarrow \C$ we have } ${\displaystyle \sum_{{d\in (Q)}\atop{d|a}}F'(d)=\sum_{{t\in (Q)}\atop{{t|a}\atop{{a\over t}\,\in\,)Q(}}}F(t) }$, $\forall a\in \N$.  
\smallskip
\par
\noindent {\bf Proof.} From the definition of Eratosthenes transform, 
$$
F'(d)=\sum_{t|d}F(t)\mu\left({d\over t}\right)
\enspace \Rightarrow \enspace
\sum_{{d\in (Q)}\atop{d|a}}F'(d)=\sum_{{t\in (Q)}\atop{t|a}}F(t)\sum_{{K\in (Q)}\atop{K\left|{a\over t}\right.}}\mu(K). 
$$ 
\par
\noindent
The characteristic function of $(Q)$ is $\1_{(Q)}$, multiplicative, so the thesis comes from 
$$
\sum_{{K\in (Q)}\atop{K|n}}\mu(K)=\sum_{K|n}\mu(K)\1_{(Q)}(K)
=\prod_{p|n}(1-\1_{(Q)}(p))
=\1_{)Q(}(n),
\enspace 
\forall n\in \N.
$$ 
\hfill $\square$ 

\medskip

\par
Our next Lemma bounds the $n\in \N$ that are $Q-$smooth and counts those which are $Q-$sifted. 
\smallskip
\par
\noindent
{\bf Lemma 2.}
\enspace {\it As } $x\to \infty$, \quad ${\displaystyle \sum_{{n\in (Q)}\atop{n\le x}}1\ll_{Q,\varepsilon}x^{\varepsilon} }$ \qquad {\it and} \qquad ${\displaystyle \sum_{{n\in )Q(}\atop{n\le x}}1=\prod_{p\le Q}\left(1-{1\over p}\right)x+O_{Q}(1) }$.  
\smallskip
\par
\noindent {\bf Proof.} We may represent (in a unique way) any $n\in(Q)$ as $n=p_1^{K_1}\cdots p_r^{K_r}$, where $2=p_1<p_2<\cdots<p_r$ are consecutive prime numbers, $K_j\ge 0$ are integers $\forall j\le r$ and this $r$ is $\pi(Q)$ (number of $p\le Q$). Then, by \lq \lq Rankin's trick\rq \rq, $\forall \varepsilon>0$ we have 
$$
\sum_{{n\in (Q)}\atop{n\le x}}1\le \sum_{{n\in (Q)}\atop{n\le x}}{{x^{\varepsilon}}\over {n^{\varepsilon}}}
\ll x^{\varepsilon} \sum_{K_1=0}^{\infty} \cdots \sum_{K_r=0}^{\infty}(p_1^{-\varepsilon})^{K_1} \cdots (p_r^{-\varepsilon})^{K_r}
=x^{\varepsilon} \prod_{p\le Q}{1\over {1-p^{-\varepsilon}}}
\ll_{Q,\varepsilon}x^{\varepsilon}. 
$$ 
\par
\noindent
This proves the bound. 
\par
\noindent
On the other side, abbreviating $P_Q:=\prod_{p\le Q}p$, the condition \enspace $(n,P_Q)=1$ \enspace is detected by \enspace ${\displaystyle \sum_{d|n,d|P_Q}\mu(d) }$: 
$$
\sum_{{n\in )Q(}\atop{n\le x}}1=\sum_{{n\le x}\atop{(n,P_Q)=1}}1
=\sum_{d\left|P_Q\right.}\mu(d)\left[{x\over d}\right]
=\sum_{d\left|P_Q\right.}{{\mu(d)}\over d}\cdot x+O\left(\sum_{d\left|P_Q\right.}\mu^2(d)\right)
=\prod_{p\le Q}\left(1-{1\over p}\right)x+O_{Q}(1). 
$$ 
\hfill $\square$ 

\medskip

\par
Our last Lemma, the core of our arguments, gives an estimate for a series restricted to $Q-$smooth numbers (badly diverging, without restrictions), that we'll use many times. (As usual,we assume $\varepsilon>0$.) 
\smallskip
\par
\noindent
{\bf Lemma 3.}
\enspace {\it For all } $0<\varepsilon<1$ {\it we get} 
$$
\sum_{m\in (Q)}m^{\varepsilon-1}\ll_{Q,\varepsilon} 1. 
$$ 
\smallskip
\par
\noindent {\bf Proof.} Representing as above the $m\in(Q)$, 
$$
\sum_{m\in (Q)}m^{\varepsilon-1}=\sum_{K_1=0}^{\infty} \cdots \sum_{K_r=0}^{\infty}(p_1^{\varepsilon-1})^{K_1} \cdots (p_r^{\varepsilon-1})^{K_r}
=\prod_{p\le Q}{1\over {1-p^{\varepsilon-1}}}
\ll_{Q,\varepsilon}1. 
$$ 
\hfill $\square$ 

\vfill
\eject

\par				
\noindent {\bf 3. Smooth-Twisted Orthogonality.}
\smallskip
\par
\noindent
We give a kind of orthogonality relations (for Ramanujan sums) which are, so to speak, smooth-twisted, i.e., they contain a kind of twist, namely the indicator function of smooth numbers (with a factor at the denominator); see that the two variables expressing the orthogonality have both to live in smooth numbers. In fact, with this restriction, the LHS (left hand side) in next result is meaningful. 
\par
We state and prove the \lq \lq Smooth-Twisted Orthogonality\rq \rq. It provides another approach to Theorem 1. 
\smallskip
\par
\noindent
{\bf Proposition 2.}
\enspace {\it Let } $q,\ell\in (Q)$. {\it Then} 
$$
{1\over {\displaystyle \sum_{t\in(Q)}{1\over t} } }
\thinspace \cdot \thinspace 
\sum_{t\in(Q)}{{c_{q}(t)c_{\ell}(t)}\over t}
=\varphi(\ell)\1_{q=\ell}\enspace. 
$$ 
\smallskip
\par
\noindent {\bf Remark 3.} We ask diligent readers to prove the absolute convergence in LHS with Ramanujan sums. 
\smallskip
\par
\noindent {\bf Proof.} Representing the denominator in LHS as 
$$
\sum_{t\in(Q)}{1\over t}=\sum_{K_1=0}^{\infty} \cdots \sum_{K_r=0}^{\infty}(p_1^{-1})^{K_1} \cdots (p_r^{-1})^{K_r}
=\prod_{p\le Q}\left(1-{1\over p}\right)^{-1}
=\left(\prod_{p\le Q}\left(1-{1\over p}\right)\right)^{-1}, 
$$ 
\par
\noindent
from the representation of numbers $t\in (Q)$, compare Lemma 2 proof, we are left with proving 
$$
\prod_{p\le Q}\left(1-{1\over p}\right)\sum_{t\in(Q)}{{c_{q}(t)c_{\ell}(t)}\over t}
=\1_{q=\ell}\varphi(\ell). 
\leqno{(\ast\ast)}
$$ 
\par
\noindent
This is a straight task, applying elementary properties like (see [M] and [T]) 
$$
c_{q}(t)=\sum_{q'|q,q'|t}q'\mu(q/q') 
\quad \hbox{\rm and} \quad 
\sum_{d|n}\varphi(d)=n, 
$$
\par
\noindent
with $n=(\ell',q')$, in the following, so to get 
$$
\prod_{p\le Q}\left(1-{1\over p}\right)\sum_{t\in(Q)}{{c_{q}(t)c_{\ell}(t)}\over t}
=\prod_{p\le Q}\left(1-{1\over p}\right)
\sum_{q'|q}\mu\left({q\over {q'}}\right)\sum_{K\in(Q)}{{c_{\ell}(q'K)}\over K}
=
$$
$$
=\prod_{p\le Q}\left(1-{1\over p}\right)
\sum_{q'|q}\mu\left({q\over {q'}}\right)
\sum_{\ell'|\ell}\ell'\mu\left({{\ell}\over {\ell'}}\right)
\sum_{{K\in(Q)}\atop{q'K\equiv 0\bmod \ell'}}{1\over K}
=
$$
$$
=\prod_{p\le Q}\left(1-{1\over p}\right)
\sum_{q'|q}\mu\left({q\over {q'}}\right)
\sum_{g'|q'}
\sum_{{\ell'|\ell}\atop{(q',\ell')=g'}}\ell'\mu\left({{\ell}\over {\ell'}}\right)
\sum_{K'\in(Q)}{1\over {K'\cdot {{\ell'}\over {g'}}}}
=
$$
$$
=
\sum_{q'|q}\mu\left({q\over {q'}}\right)
\sum_{g'|q'}g'
\sum_{{\ell'|\ell}\atop{(q',\ell')=g'}}\mu\left({{\ell}\over {\ell'}}\right)
=\sum_{q'|q}\mu\left({q\over {q'}}\right)
\sum_{\ell'|\ell}\mu\left({{\ell}\over {\ell'}}\right)(\ell',q')
=\sum_{q'|q}\mu\left({q\over {q'}}\right)
\sum_{\ell'|\ell}\mu\left({{\ell}\over {\ell'}}\right)\sum_{{d|\ell'}\atop{d|q'}}\varphi(d)
=
$$
$$
=
\sum_{q'|q}\mu\left({q\over {q'}}\right)
\sum_{{d|\ell}\atop{d|q'}}\varphi(d)\sum_{{\ell'|\ell}\atop{\ell'\equiv 0\bmod d}}\mu\left({{\ell}\over {\ell'}}\right)
=\sum_{q'|q}\mu\left({q\over {q'}}\right)
\sum_{{d|\ell}\atop{d|q'}}\varphi(d)\sum_{\ell''|{{\ell}\over d}}\mu\left({{\ell/d}\over {\ell''}}\right)
=\sum_{q'|q}\mu\left({q\over {q'}}\right)\varphi(\ell)\1_{\ell|q'}
=
$$
$$
=
\varphi(\ell)\1_{\ell|q}\sum_{{q'|q}\atop{q'\equiv 0\bmod \ell}}\mu\left({q\over {q'}}\right)
=\varphi(\ell)\1_{\ell|q}\sum_{q''\left|{q\over {\ell}}\right.}\mu\left({{q/\ell}\over {q''}}\right)
=\1_{q=\ell}\varphi(\ell), 
$$
\par
\noindent
since [T] M\"obius inversion \enspace ${\displaystyle \sum_{d|n}\mu(d) }=\1_{\{1\}}(n)$ \enspace is applied twice. 
Thus $(\ast\ast)$ is completely proved.\hfill $\square$ 
 
\vfill
\eject

\par				
\noindent {\bf 4. A simple counterexample to the Reef.}
\smallskip
\par
\noindent
Previous version 2 of present paper contains the two following Conjectures. 
\smallskip
\par
\noindent
{\bf Conjecture 1.}
\enspace {\it Let } $q,\ell\in (Q)$ {\it and } $n\in \Z$. {\it Then} 
$$
{1\over {\displaystyle \sum_{t\in(Q)}{1\over t} } }
\thinspace  \thinspace 
\sum_{t\in(Q)}{{c_{q}(n+t)c_{\ell}(t)}\over t}
\enspace=\enspace\1_{q=\ell}c_{\ell}(n). 
$$ 
\smallskip
\par
\noindent {\bf Remark 4.} If we take any $n\equiv 0(\!\bmod \; q)$ we get previous Proposition 2. 
\smallskip
\par
\noindent
{\bf Conjecture 2.}
\enspace {\it Let } $f,g$ {\it satisfy } $(BH)$. {\it Then we have for their correlation the } $Q-${\it smooth restricted Reef} 
$$
C_{f,g}(N,a)=\sum_{\ell\le Q}\left({{\widehat{g}(\ell)}\over {\varphi(\ell)}}\sum_{n\le N}f(n)c_{\ell}(n)\right)c_{\ell}(a), 
\enspace \forall a\in (Q).
\leqno{(Q)-\hbox{\stampatello Reef}:}
$$ 
\smallskip
\par 
In version 2 we prove that Conjecture 1 implies Conjecture 2, i.e. the following. 
\smallskip
\par
\noindent
{\bf Proposition 3.}
\enspace {\it Let } $f,g$ {\it satisfy } $(BH)$ {\it and let Conjecture 1 hold. Then we have the } $(Q)-${\stampatello Reef}. 

\bigskip

We give an important counterexample to the Reef, namely we disprove now Conjecture 2 (hence, thanks to Proposition 3, disproving Conjecture 1, too). 
\smallskip
\par
\noindent
{\bf Counterexample 1.}
\enspace {\it Let } $N,Q\in \N$ {\it and fix the two integers } $1\le n_0\le N$ {\it and } $2<q_0\le Q$. {\it Then, choosing}
$$
f(n)\defineq \1_{\{n_0\}}(n),
\enspace \forall n\in \N 
\qquad 
\hbox{\it and}
\qquad
g(m)\defineq c_{q_0}(m),
\enspace \forall m\in \N 
$$
\par
\noindent
{\it we have $(BH)$ for $f$ and $g$, but we can't have the } $(Q)-${\stampatello Reef} {\it since}
$$
a=1, n_0\equiv-1(\bmod \thinspace \thinspace q_0)
\enspace \Rightarrow \enspace 
C_{f,g}(N,a)=\varphi(q_0)
\neq
{1\over {\varphi(q_0)}}\mu^2(q_0)=\sum_{\ell\le Q}\left({{\widehat{g}(\ell)}\over {\varphi(\ell)}}\sum_{n\le N}f(n)c_{\ell}(n)\right)c_{\ell}(a). 
$$

\bigskip

\par
\noindent
Since, for general $n_0,q_0$ (in the above hypotheses), Counterexample 1 has
$$
C_{f,g}(N,a)=c_{q_0}(n_0+a)
\qquad 
\hbox{\rm and}
\qquad
\sum_{\ell\le Q}\left({{\widehat{g}(\ell)}\over {\varphi(\ell)}}\sum_{n\le N}f(n)c_{\ell}(n)\right)c_{\ell}(a)
={1\over {\varphi(q_0)}}c_{q_0}(n_0)c_{q_0}(a)
$$
\par
\noindent
for all $a\in \N$, then, at least for this case, we may substitute Conjecture 2 with : 
\smallskip
\par
\noindent
{\bf Conjecture 3.}
\enspace {\it Let } $f,g$ {\it satisfy } $(BH)$ {\it and assume $\exists \delta>0$ such that } $Q\le N^{1-\delta}$. {\it Then we have the} 
$$
C_{f,g}(N,a)=\sum_{\ell\le Q}\left({{\widehat{g}(\ell)}\over {\varphi(\ell)}}\sum_{n\le N}f(n)c_{\ell}(n)\right)c_{\ell}(a)
             +O_{\delta}\left(N^{1-\delta}\right), 
\enspace \forall a\le N^{1-\delta}.
\leqno{\hbox{\stampatello Approximate\enspace Reef}:}
$$ 
\smallskip
\par
\noindent {\bf Remark 5.} We are assuming a very strong remainder and, also, a very large range of uniformity for $a$. 

\vfill
\eject

\par				
\noindent {\bf 5. Further remarks.}
\smallskip
\par
\noindent
Since $V\ge Q$ and $\ell \le Q$ imply $\ell\in(V)$, from our two Corollaries above we easily get our next result. 
\smallskip
\par
\noindent {\bf Corollary 3.}
\enspace {\it Given } $f:\N \rightarrow \C$ {\it and } $g:\N \rightarrow \C$ {\it satisfying the Basic Hypothesis $(BH)$ above, we have} 
$$
V\ge Q 
\enspace 
\Rightarrow 
\quad
\forall a\in (V),
\enspace 
C_{f,g}(a)=\sum_{\ell\le Q}{{\widehat{g}(\ell)}\over {\varphi(\ell)}}\sum_{n\le N}f(n)c_{\ell}(n)c_{\ell}(a)
           -\sum_{\ell\in(V)}\sum_{{d\not\in(V)}\atop{d\equiv 0(\!\!\bmod \ell)}}{{C'_{f,g}(N,d)}\over d}c_{\ell}(a),
$$
\par
\noindent
{\it whence, passing to the limit over } $V\in \N$, {\it we obtain } $\forall a\in \N$ \enspace {\it the}
$$
C_{f,g}(N,a)=\sum_{\ell\le Q}{{\widehat{g}(\ell)}\over {\varphi(\ell)}}\sum_{n\le N}f(n)c_{\ell}(n)c_{\ell}(a)
             -\lim_{V}\sum_{\ell\in(V)}\sum_{{d\not\in(V)}\atop{d\equiv 0(\!\!\bmod \ell)}}{{C'_{f,g}(N,d)}\over d}c_{\ell}(a). 
\leqno{\hbox{\stampatello Asymptotic\enspace Reef}:}
$$ 

\bigskip

\par
\noindent
The present results have, of course, applications to our study in [C1], [C2], [CL] and in the series of papers starting with [CMS], [CM]. In particular, they may be applied to averages of correlations (see [CL]) and to single correlations [C2], [CM], with a more expected success (for reasons that we'll explain in future papers) for the averages (having, see [CL], a big impact on moments of the Riemann $\zeta-$function on the critical line). 

\medskip

\par
For a more extensive discussion on these arguments, compare especially Generations [CL] and [CM]. For remarks on the Ramanujan expansion coefficients and their decay see [C2] and [CM]. Last but not least, for applications to conditional proofs of Hardy-Littlewood Conjecture, compare [C1]. 

\medskip

\par
\noindent
{\bf Acknowledgments.} I wish to thank again Ram Murty for our previous common papers, a source of inspiration for [C1], [C2] and present paper. Also, I wish to thank Maurizio Laporta for an extensive and careful reading of previous versions of the paper. 

\bigskip

\par
\centerline{\stampatello Bibliography}

\smallskip

\item{[Ca]} R.D. Carmichael, {\sl Expansions of arithmetical functions in infinite series}, Proc. London Math. Society {\bf 34} (1932), 1--26. 


\item{[C1]} G. Coppola, {\sl An elementary property of correlations}, Hardy-Ramanujan J. {\bf 41} (2018), 68--76. Available online 


\item{[C2]} G. Coppola, {\sl A map of Ramanujan expansions}, ArXiV:1712.02970v2. (Second Version) 


\item{[CL]} G. Coppola and M. Laporta, {\sl Generations of correlation averages}, J. Numbers Volume 2014 (2014), Article ID 140840, 13 pages http://dx.doi.org/10.1155/2014/140840 (draft, ArXiV:1205.1706.) 


\item{[CMS]} G. Coppola, M. Ram Murty and B. Saha, {\sl Finite Ramanujan expansions and shifted convolution sums of arithmetical functions},  J. Number Theory {\bf 174} (2017), 78--92. 


\item{[CM]} G. Coppola and M. Ram Murty, {\sl Finite Ramanujan expansions and shifted convolution sums of arithmetical functions, II}, J. Number Theory {\bf 185} (2018), 16--47. 


\item{[Da]} H. Davenport, {\sl Multiplicative Number Theory}, 3rd ed., GTM 74, Springer, New York, 2000. 


\item{[De]} H. Delange, {\sl On Ramanujan expansions of certain arithmetical functions}, Acta Arith. {\bf 31}(1976), 259--270. Available online 


\item{[De87]} H. Delange, {\sl On a formula for almost-even arithmetical functions}, Illinois J. Math. {\bf 31} (1987), 24--35. Available online 


\item{[M]} M. Ram Murty, {\sl Ramanujan series for arithmetical functions}, Hardy-Ramanujan J. {\bf 36} (2013), 21--33. Available online 


\item{[R]} S. Ramanujan, {\sl On certain trigonometrical sums and their application to the theory of numbers}, Transactions Cambr. Phil. Soc. {\bf 22} (1918), 259--276. 


\item{[T]} G. Tenenbaum, {\sl Introduction to Analytic and Probabilistic Number Theory}, Cambridge Studies in Advanced Mathematics, {46}, Cambridge University Press, 1995. 


\item{[W]} A. Wintner, {\sl Eratosthenian averages}, Waverly Press, Baltimore, MD, 1943. 

\bigskip

\leftline{\tt Giovanni Coppola - Universit\`{a} degli Studi di Salerno (affiliation)}
\leftline{\tt Home address : Via Partenio 12 - 83100, Avellino (AV) - ITALY}
\leftline{\tt e-mail : giovanni.coppola@unina.it}
\leftline{\tt e-page : www.giovannicoppola.name}
\leftline{\tt e-site : www.researchgate.net}

\vfill
\eject

\bye